\documentclass[letterpaper,12pt]{amsart}
\title[Decompositions of Betti diagrams of powers of monomial ideals]{Decompositions of Betti diagrams of powers of monomial ideals: A stability conjecture
\\ $\,$ \\  {\tiny Extended abstract for the conference\\ combinatorial methods in topology and algebra, Cortona 2013 }}

\author{Alexander Engstr\"om}
\address{Department of Mathematics \\ Aalto University}
\email{alexander.engstrom@aalto.fi}
\date{\today}

\newtheorem{conjecture}{Conjecture}

\begin{document}

\maketitle

\section{The conjecture}
We state a conjecture on the stability of Betti diagrams of powers of monomial ideals. Boij and S\"oderberg  \cite{AE_BS} conjectured that Betti diagrams can be decomposed into pure diagrams, and that was proved by Eisenbud and Schreyer \cite{AE_ES}. We don't cover the basics of Boij-S\"oderberg theory in this extended abstract, see Fl\o ystad \cite{AE_F} for a survey.

Up to scaling a pure diagram is determined by its non-zero positions. In our setting the top left corner is always non-zero and we normalise by assigning it the value one. For higher powers of ideals we need taller pure diagrams in a sequential way. A \emph{translation of a pure diagram} for $k=0,1,2,\ldots$ is a sequence of pure diagrams on the form
\[
\begin{array}{rcl}
\pi(k) & = &
\begin{array}{r|l}
& 0 \,\,\,\,\,  1 \,\,\,\,\,  2 \,\,\,\,\,  \cdots  \,\,\,\,\,   \\
\hline
0 & 1 \\
\vdots \\
l(k) &   \,\,\,\,\,  \,\,\,\,\, \textrm{A fixed shape for} \\ 
& \,\,\,\,\,  \,\,\,\,\, \textrm{non-zero entries} \\ 
\end{array}
\end{array}
\]
where $l(k)$ is a linear function.

According to Boij-S\"oderberg theory there is for every ideal $I$ in $S$ a decomposition of the Betti diagram $\beta(S/I)=w_1\pi_1+\cdots +w_m\pi_m$ where each $w_i$ is a non-negative real number and each $\pi_i$ is a pure diagram. Usually there are many choices of weights, and when considering algorithms to find decompositions there is a point to finding a particular one. But the amount of choices is also a measure of the complexity of the Betti diagram, and one might notice that for ideals that we know the invariants of for large powers, the complexity in this sense is quite low. For example, if all powers are linear, then there are no choices at all.

For any $\beta(S/I)$ there is a finite set of pure diagrams that can be included in a decomposition with a positive weights. We call the set of possible weight vectors for $\beta(S/I)$ the \emph{polytope of Betti diagram decompositions}. 

\newpage
We conjecture that for high powers the polytope of Betti diagram decompositions stabilises.

\begin{conjecture}
Let $I$ be a monomial ideal in $S$ with all generators of the same degree. Then there is a $k_0$ such that for all $k>k_0,$
\begin{itemize}
\item[1.] For some translations of pure diagrams $\pi_1(k),\ldots,\pi_m(k)$ any decomposition of $\beta(S/I^k)$ is a weighted sum on the form  $w_1\pi_1(k)+\cdots+w_m\pi_m(k).$ Denote the polytope of Betti diagram decompositions of $\beta(S/I^k)$ in $\mathbb{R}^m$ by $P_k.$
\item[2.] All $P_k$ are of the same combinatorial type as a polytope $P_I$. For any vertex $v$ of $P_I$ there is a function $h_v(k)\in \mathbb{R}^m$, which is rational in each coordinate, such that the vertex corresponding to $v$ in $P_k$ is $h_v(k).$
\end{itemize}
\end{conjecture} 

The conjecture is true for ideals whose large enough powers are all linear: The polytope is a point. This follows from that the column sums in Betti diagrams stabilises to polynomials for large powers according to Kodiyalam \cite{AE_K} and from the procedure to derive the unique decomposition of linear diagrams. The conjecture holds for many small examples that the author have calculated. There is unfortunately no abundance of ideals in the literature for which the Betti diagrams of all powers are given explicitly, since these concepts are fairly new. But there are many interesting tools accessible, for example from algebraic and topological combinatorics, that should make serious attempts to derive them fruitful.

\section{An example}

In this section we give an example of an ideal satisfying the conjecture. 
Engstr\"om and Nor\'en \cite{AE_EN} constructed explicit cellular minimal resolutions of $S/I^k$ for all $k$ and $n,$ where
\[ S=\mathbf{k}[x_1,x_2,\ldots,x_n] \,\,\, \textrm{and} \,\,\, I=\langle x_1x_2,x_2x_3,\ldots,x_{n-1}x_n \rangle, \]
and calculated the Betti numbers:
\[
\beta_{i,j}(S/I^k)={n+3k-j-2 \choose 2j-3i-3k+3}{n+4k+2i-2j-4 \choose 2k+2i-j-2}{j-i-k \choose k-1}.
\]

The Betti diagram of $\mathbf{k}[x_1,x_2,x_3,x_4,x_5,x_6] / \langle x_1x_2,x_2x_3,x_3x_4,x_4x_5,x_5x_6 \rangle^k$ is
\[
\begin{array}{c|cccccc}
& 0 & 1 & 2 & 3 & 4 & 5 \\ 
\hline
0 & 1 \\
\vdots  \\
2k-1 && {k+4 \choose 4} & 4{k+3 \choose 4} & 6{k+2 \choose 4} & 4{k+1 \choose 4} & {k \choose 4} \\
&&& k(k+2) & 2k(k+1) & k^2 & \\
\end{array}
\]

\newpage
The translations of pure diagrams:

\[\begin{array}{cc}
\begin{array}{rcl}
\pi_1(k) & = &
\begin{array}{c|cccc}
& 0 & 1 & 2 & 3 \\
\hline
0 & 1 \\
\vdots \\
2k-1 & & \ast & \ast & \ast \\
\end{array}
\end{array}
&
\begin{array}{rcl}
\pi_2(k) & = &
\begin{array}{r|cccc}
& 0 & 1 & 2 & 3 \\
\hline
0 & 1 \\
\vdots \\
2k-1 & & \ast & \ast &  \\
&&&& \ast \\
\end{array}
\end{array} \\
\, \\
\begin{array}{rcl}
\pi_3(k) & = &
\begin{array}{c|cccc}
& 0 & 1 & 2 & 3 \\
\hline
0 & 1 \\
\vdots \\
2k-1 & & \ast &&  \\
&&& \ast & \ast \\
\end{array}
\end{array} &
\begin{array}{rcl}
\pi_4(k) & = &
\begin{array}{c|ccccc}
& 0 & 1 & 2 & 3 & 4 \\
\hline
0 & 1 \\
\vdots \\
2k-1 & & \ast & \ast & \ast & \ast\\
\end{array}
\end{array} \\
\, \\
\begin{array}{rcl}
\pi_5(k) & = &
\begin{array}{r|ccccc}
& 0 & 1 & 2 & 3 & 4 \\
\hline
0 & 1 \\
\vdots \\
2k-1 & & \ast & \ast & \ast \\
&&&&& \ast
\end{array}
\end{array} & 
\begin{array}{rcl}
\pi_6(k) & = &
\begin{array}{c|ccccc}
& 0 & 1 & 2 & 3 & 4 \\
\hline
0 & 1 \\
\vdots \\
2k-1 & & \ast & \ast & \\
&&&& \ast & \ast
\end{array}
\end{array} \\
\, \\
\begin{array}{rcl}
\pi_7(k) & = &
\begin{array}{r|ccccc}
& 0 & 1 & 2 & 3 & 4 \\
\hline
0 & 1 \\
\vdots \\
2k-1 & & \ast &  \\
&&& \ast & \ast & \ast
\end{array}
\end{array} &
\begin{array}{rcl}
\pi_8(k) & = &
\begin{array}{r|cccccc}
& 0 & 1 & 2 & 3 & 4 & 5\\
\hline
0 & 1 \\
\vdots \\
2k-1 & & \ast & \ast & \ast & \ast & \ast \\
\end{array}
\end{array}
\end{array}\]

The polytope $P_k$ is a triangle whose vertices have the coordinates $h_1(k),$ $h_2(k)$ and $h_3(k).$

\[h_1(k)=\left(
0,0, \frac{k+2}{2k+3}, w_4(k),
\frac{(2k+5)(k-1)}{(2k+1)(k+2)(k+1)}, \frac{(4k+5)(k+1)}{(2k+3)(2k+1)(k+2)},0,w_8(k)
\right) \]
\[h_2(k)=\left(
0,\frac{2(k+2)(k+2)}{(2k+3)(2k+1)},0, w_4(k),
\frac{(2k+5)(k-1)}{(2k+1)(k+2)(k+1)}, \frac{(k+1)(k-1)}{(2k+3)(2k+1)(k+2)},\frac{1}{2k+3},w_8(k)
\right) \]
\[h_3(k)=\left(
\frac{k+2}{2k+1},0,0, w_4(k),
\frac{k^2-7}{(2k+1)(k+2)(k+1)}, \frac{k+1}{(2k+3)(2k+1)(k+2)},\frac{1}{2k+3},w_8(k)
\right) \]
\[w_4(k)=  \frac{(7k+5)(k-1)(k-2)}{4(2k+3)(2k+1)(k+1)} \]
\[w_8(k)= \frac{(k-1)(k-2)(k-3)}{4(2k+3)(2k+1)(k+1)} \]

\end{document}